\documentclass[sn-mathphys,Numbered]{sn-jnl}

\usepackage{graphicx}%
\usepackage{multirow}%
\usepackage{amsmath,amssymb,amsfonts}%
\usepackage{amsthm}%
\usepackage{mathrsfs}%
\usepackage[title]{appendix}%
\usepackage{xcolor}%
\usepackage{textcomp}%
\usepackage{manyfoot}%
\usepackage{booktabs}%
\usepackage{algorithm}%
\usepackage{algorithmicx}%
\usepackage{algpseudocode}%
\usepackage{listings}%
\usepackage{comment}%
\usepackage[justification=centering]{caption}%

\theoremstyle{plain}%
\theoremstyle{definition}%
\newtheorem{conjecture}{Conjecture}%

\theoremstyle{remark}%

\begin{document}

\title[$K_5$ and $K_{3,3}$ are Toroidal Penny Graphs]{$K_5$ and $K_{3,3}$ are Toroidal Penny Graphs}


\author*[1]{ \fnm{Cédric} \sur{Lorand}}\email{cedric.lorand@gmail.com}


\abstract{In this article, we emphasize the connection between two fields of study, namely penny graphs and the optimal packing of spheres on the flat torus. We give a brief literature overview of related results in planar graphs, penny graphs, toroidal penny graphs, and spherical codes. We also show that $K5$ and $K_{3,3}$ are penny graphs on the flat square torus.}


\maketitle

\section{Introduction}\label{sec1}

The study of planar graphs has evolved significantly over the past century, starting with Kuratowski’s theorem, which ensures planarity by excluding  subdivisions of $K_5$ and $K_{3,3}$. Since the 1970's, the development of linear time planarity testing algorithms further deepened our understanding.

Penny graphs, formed by tangent unit circles, gained interest in the 1970s, leading to findings on edge bounds and the NP-hardness of embedding arbitrary graphs as penny graphs. Following earlier works on spherical codes, recent research on optimal sphere packings of the flat torus emerged. While advancements have been made in understanding optimal circle packings, proving their optimality continues to pose significant challenges.

\section{\texorpdfstring{$K_5$, $K_{3,3}$}{Lg} and Planar Graphs}\label{sec2}

\begin{figure}[ht]%
\centering
\includegraphics[width=0.9\textwidth]{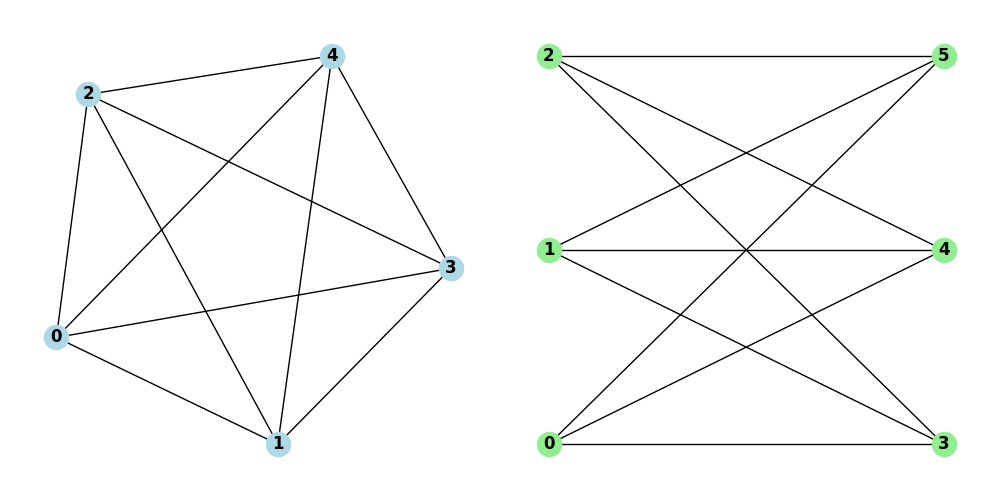}
\caption{ \textbf{left}: the complete graph $K_5$, \textbf{right}: the complete bipartite graph $K_{3,3}$}
\label{fig1}
\end{figure}

For over a century, planar graphs have been extensively studied. In 1930, Kuratowski published a proof of his now-famous theorem.\cite{bib0} It states that an arbitrary graph is planar if and only if, it does not contain as a subgraph any subdivision of either: $K_5$ the complete graph with 5 vertices; or $K_{3,3}$  the complete bipartite graph with 6 vertices (cf. Figure~\ref{fig1}). 

Later, in 1974, Hopcroft and Tarjan proposed a powerful linear time algorithm that detects planarity of a graph in $O\left(V\right)$, 
where $V$ is the number of vertices in the graph. \cite{bib1} The algorithm does not rely on Kuratowski's theorem. Instead, it uses an iterative depth-first search algorithm as follows.
First, it extracts an initial cycle from the input graph. Then, it removes the cycle from the graph which creates a set of connected components. Finally, it attempts
to iteratively add the components to the cycle, in such a way that the obtained graph remains planar. If it at a certain point it is impossible to add a component
and keep planarity, the algorithm returns "non-planar", otherwise after all components were added it returns "planar".  

This result is a rather surprising one. Indeed,  we know today that the problem of extracting a particular minor subgraph from an 
arbitrary graph is NP-hard in general. However,  surprisingly,  for the particular case of $K_5$ and $K_{3,3}$, planarity testing  (in linear time)
is sufficient to identify that they are not included in the input graph.

\section{Planar Penny Graphs}\label{sec3}

\begin{figure}[h]%
\centering
\includegraphics[width=0.5\textwidth]{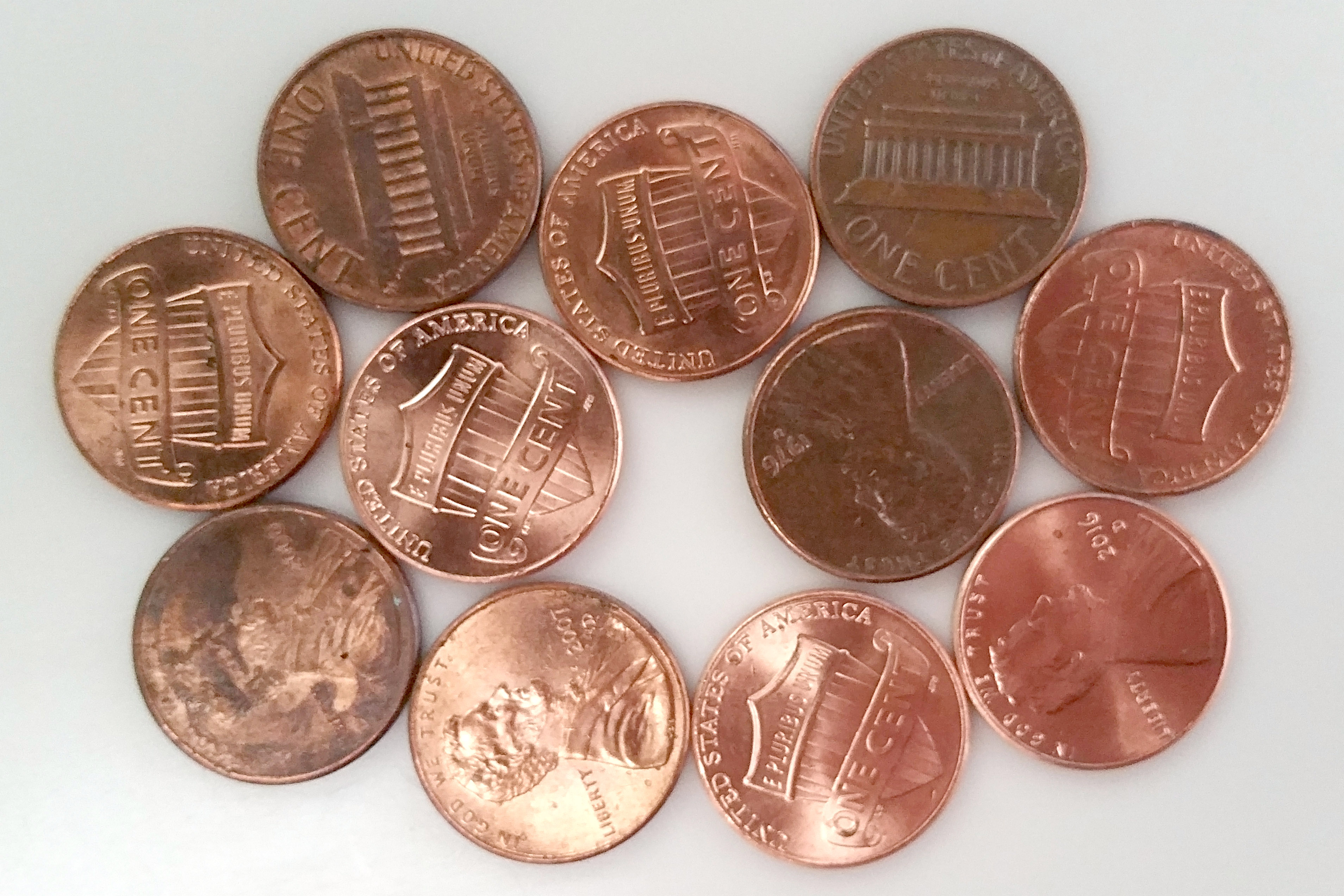}
\caption{\href{https://commons.wikimedia.org/w/index.php?curid=56426404}{11 pennies, forming a penny graph with 11 vertices and 19 edges. By David Eppstein - Own work, CC BY-SA 4.0} }
\label{fig2}
\end{figure}

Planar penny graphs are defined as the contact graph formed by tangent unit circles on the plane  (cf. Figure~\ref{fig2}). As a field of research interest has been limited. Nevertheless, since the 1970s, some interesting results have been discovered.

In 1974, Harborth proved that a penny graph has at most $\lfloor 3n-\sqrt{12n-3} \rfloor$ edges. \cite{bib2} The $3n$ term in this formula is related to the hexagonal lattice, in which there are 
exactly 3 edges per vertex. The $\sqrt{12n-3}$ correcting term accounts for the fact that on the graph convex hull boundary, vertices have less neighbors. This bound can be reached exactly 
on penny graphs that are subsections of the hexagonal lattice.

In 1990, Eades and Wormald showed that the problem of deciding whether an arbitrary graph can be embedded on the plane as a penny graph is NP-hard. \cite{bib3} Later on in 1996, this result
was generalized for the case of nearest neighbours graph  by Eades and Whitesides, a less restrictive type of graph. \cite{bib4}

\section{Optimal Spheres' Packings on the Flat Torus}\label{sec4}

Until recently the problem of penny graphs embeddings on the flat torus did not receive much attention in the literature.
Indeed, mathematicians did not directly address it from a purely graph-theoretical perspective. Rather, they approached it as an optimal sphere packing problem on the flat torus.
 This is probably because it was introduced as a generalization of ‘Spherical Codes’ also known as the ‘Tammes’ Problem. In the 1930s, Tammes began researching this problem.\cite{bib9}.
HE was followed by Leech \cite{bib10} and Kottwitz \cite{bib12} among others (\cite{bib11}, \cite{bib13}, \cite{bib14}).
Online resources on the subject of Spherical Code are available (cf Sloane  \& all \cite{bib15}).
It is also the focus of a chapter in the book by Conway and Sloane \cite{bib16}.

As stated previously, research on toroidal penny graphs, or equivalently optimal spheres packing on the flat torus is rather recent.
In 2011 Dickinson and al presented the solutions for optimal packings up to 5 circles. \cite{bib8} In 2016, Musin and Nikitenko presented solutions up to 8 circles, as well as a conjecture for the case of 9 circles. \cite{bib6}  In 2016, Conelly and al extended the results with a set of conjectures up to 16 circles. \cite{bib7} With the help of computers and the use of numerical algorithms it seems possible to exhibit putative solutions for a small number of circle. However, given a particular packing solution, proving that it is optimal is a hard problem, as generally many suboptimal solutions also exist.

\section{\texorpdfstring{$K_5$ and $K_{3,3}$}{Lg} as Toroidal Penny Graphs }\label{sec6}

\subsection{\texorpdfstring{$K_5$}{Lg} is a Toroidal Penny Graphs}\label{subsec61}
	Dickinson and al proved that the penny graph shown in Figure~\ref{fig3} is the optimal packing solution for 5 circles on the flat square torus. \cite{bib8}
 For the rest of this paper, we express toroidal coordinates such as in Table~\ref{tab1} with respect to the unit torus centered at the origin, since this greatly simplifies the formulas.

\begin{figure}[h]%
\centering
\includegraphics[width=0.9\textwidth]{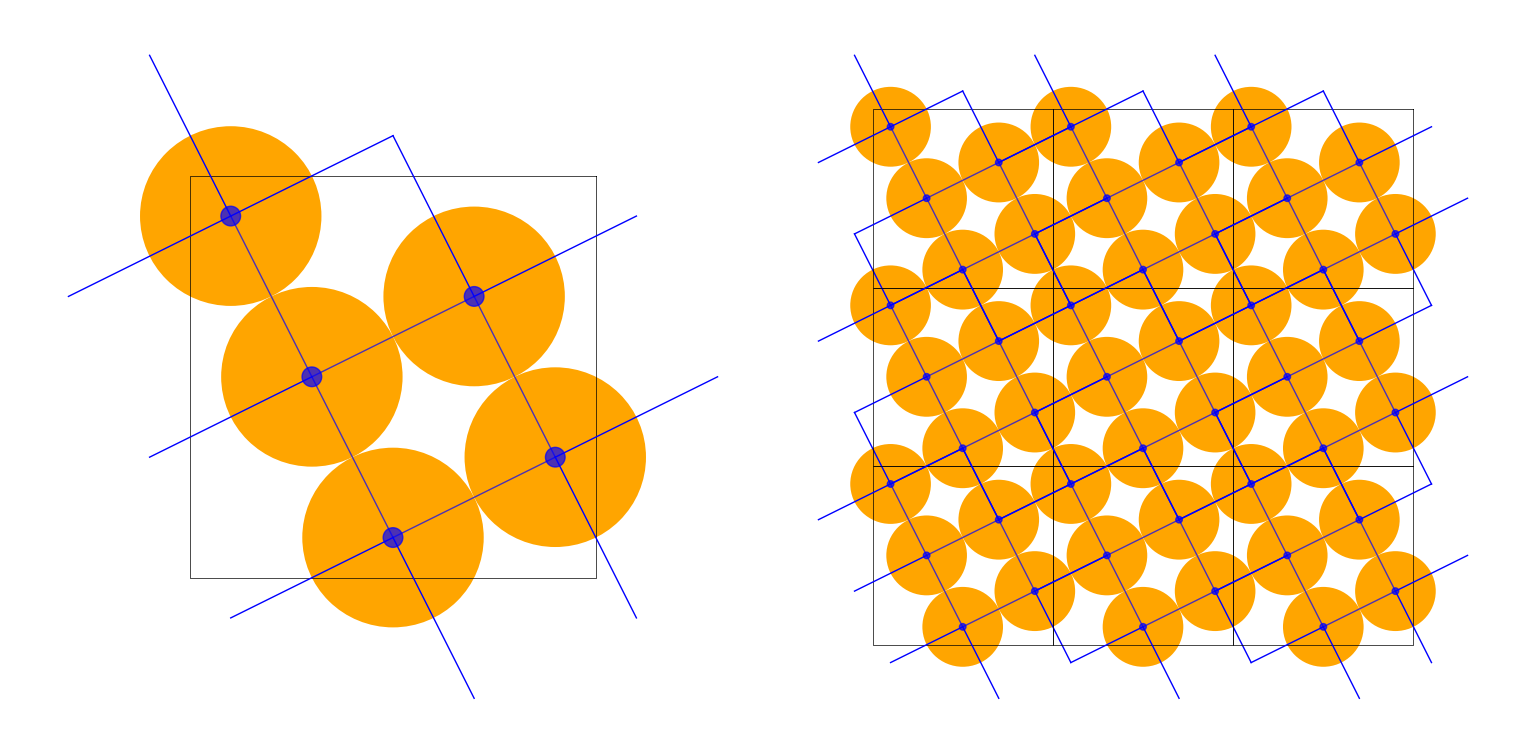}
\caption{\textbf{left}:  $K_5$ penny graph embedding on the unit flat square torus, \\ \textbf{right}: $K_5$  penny graph embedding on a 3x3 toroidal tiling}
\label{fig3}
\end{figure}

\begin{table}[ht]
\begin{tabular}{@{}llllll@{}}
\toprule
\textbf{Node \#} & \multicolumn{1}{c}{1} & \multicolumn{1}{c}{2} & \multicolumn{1}{c}{3} & \multicolumn{1}{c}{4} & \multicolumn{1}{c}{5} \\
\midrule
\textbf{$\left(x, y\right)$} & $\left(-\frac{2}{5}, \frac{2}{5}\right)$ & $\left(-\frac{1}{5}, 0\right)$ &  $\left(0, -\frac{2}{5}\right)$  & $\left(\frac{2}{5}, -\frac{1}{5}\right)$ & $\left(\frac{1}{5}, \frac{1}{5}\right)$  \\ 
\botrule
\end{tabular}
\caption{$K_5$ toroidal embedding coordinates}
\label{tab1}
\end{table}

  From the coordinates in Table~\ref{tab1} one can easily verify that the value of the packing diameter is $l=\frac{1}{\sqrt{5}}$. Moreover, note that the contact graph is the 4-regular graph with 5 nodes also known as $K_5$. This shows that $K_5$ is not only embeddable without crossing edges on the torus but it is also as a toroidal penny graph.

\subsection{Optimal Toroidal Packing of 6 Spheres}\label{subsec62}

Musin and Nikitenko showed that the packing in Figure~\ref{fig4} is the optimal packing solution for 6 circles on the flat square torus. \cite{bib6} Toroidal coordinates for this packing are given in Table~\ref{tab2}. Here $l$ stands for the packing diameter $l=\frac{1}{6}\left(1+3\sqrt{3}-\sqrt{4+6\sqrt{3}}\right)$ \cite{bib6}.

From the coordinates in Table~\ref{tab2} one can easily verify that all edges have identical length equal to the packing diameter $l$.
Here, the contact graph is the 4-regular graph with 6 nodes also known as the octahedral graph. Interestingly, the octahedral graph happens to be a planar graph, while the optimal contact graph for 5 nodes $K_5$ is not.

\begin{figure}[ht]%
\centering
\includegraphics[width=0.9\textwidth]{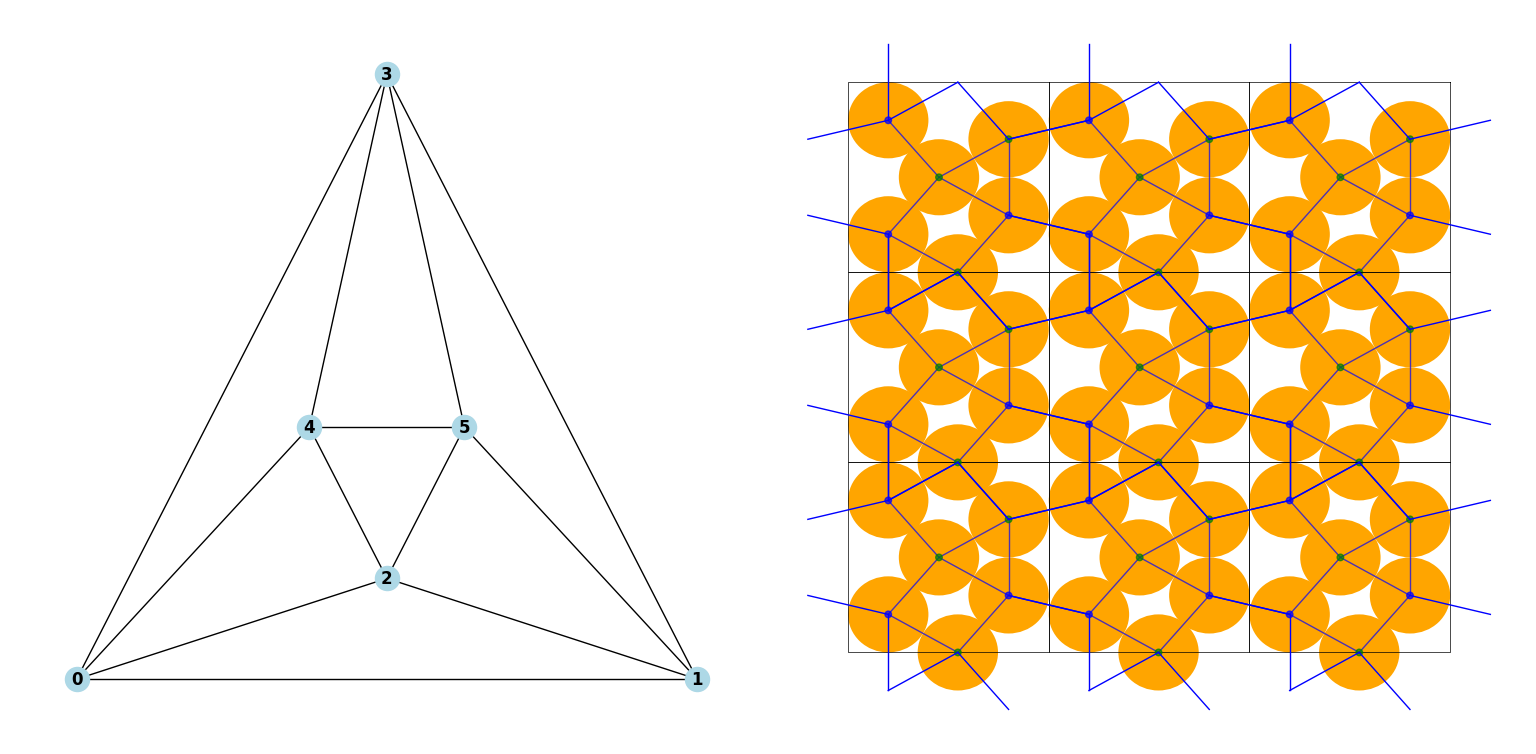}
\caption{\textbf{left}:  Planar embedding of the octahedral graph, \\ \textbf{right}: Penny graph embedding of the octahedral graph}\label{fig3}
\end{figure}

\begin{table}[ht]
\begin{tabular}{@{}lllllll@{}}
\toprule
\textbf{Node \#} & \multicolumn{1}{c}{1} & \multicolumn{1}{c}{2} & \multicolumn{1}{c}{3} & \multicolumn{1}{c}{4} & \multicolumn{1}{c}{5} & \multicolumn{1}{c}{6}\\
\midrule
\textbf{$\left(x, y\right)$} & 
$\left(\frac{l-1}{2}, \frac{l-1}{2}\right)$\footnotemark & 
$\left(\frac{1+\sqrt{3}}{2}l-\frac{1}{2}, -\frac{1}{2}\right)$ &  
$\left(-\frac{l-1}{2},-\frac{l}{2}\right)$  & 
$\left(-\frac{l-1}{2}, \frac{l}{2}\right)$  & 
$\left( \frac{l-1}{2},-\frac{l-1}{2}\right)$ & 
$\left(-\frac{1+\sqrt{3}}{2}l+\frac{1}{2},0\right)$\\ 
\botrule
\end{tabular}
\caption{Toroidal embedding coordinates for the toroidal octahedral penny graph}
\label{tab2}
\end{table}

\subsection{\texorpdfstring{$K_{3,3}$}{Lg} is a Toroidal Penny Graphs}\label{subsec7}

\begin{figure}[ht]%
\centering
\includegraphics[width=0.9\textwidth]{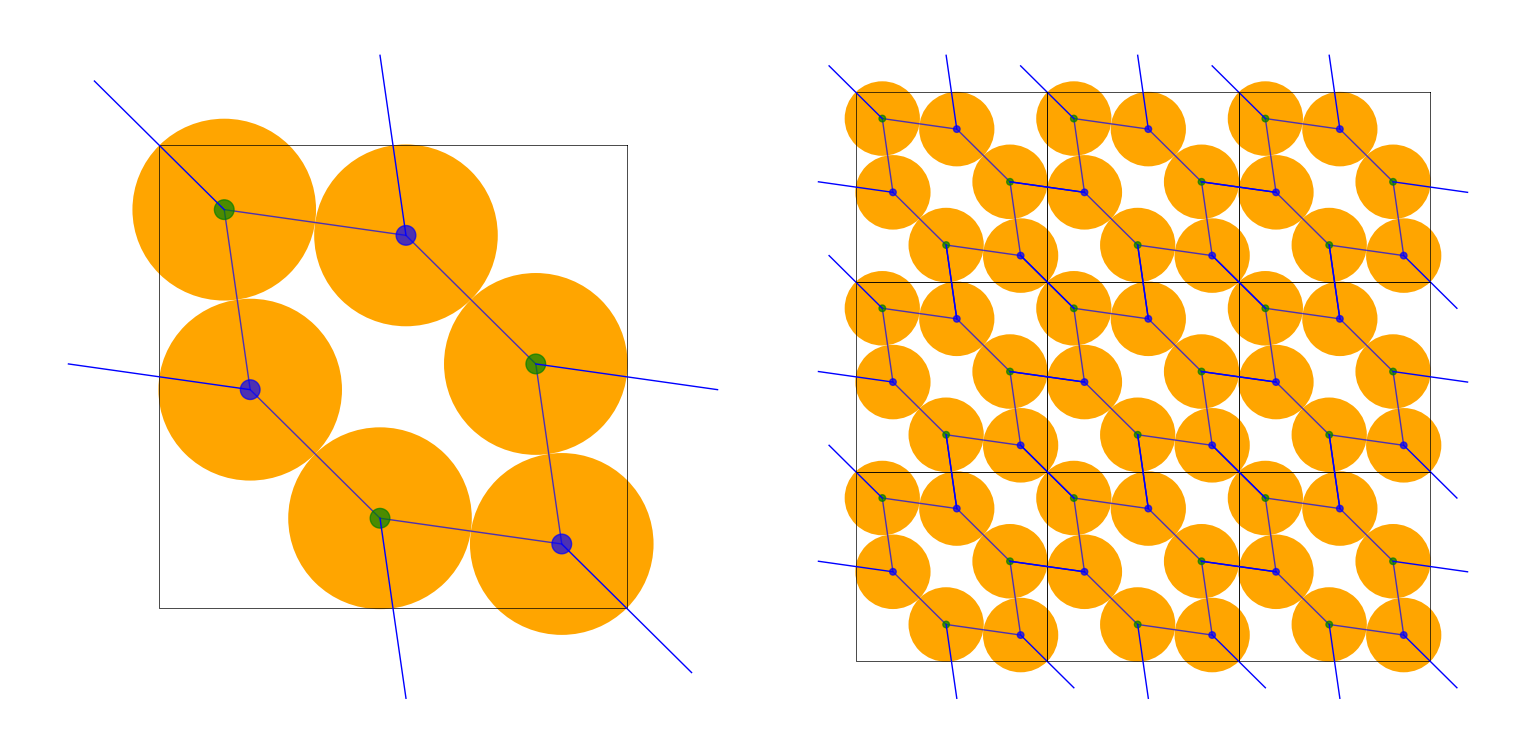}
\caption{\textbf{left}:  $K_{3,3}$ penny graph embedding on the unit flat square torus, \\ \textbf{right}: $K_{3,3}$  penny graph embedding on a 3x3 toroidal tiling}
\label{fig4}
\end{figure}

Due to its prevalent role in Kuratowski’s theorem we now investigate the question of the embedabbility of $K_{3, 3}$ as a toroidal penny graph. We claim that it is indeed a toroidal penny graph by exhibiting its embedding in Figure~\ref{fig4}. Toroidal coordinates of the nodes are provided in Table~\ref{tab3}. Once again, given the coordinates in this table one can easily verify that all edges' lengths are equal, and that the packing radius is equal to $l=\frac{5\sqrt{2}}{18}$.

\begin{table}[ht]
\begin{tabular}{@{}lllllll@{}}
\toprule
\textbf{Node \#} & \multicolumn{1}{c}{1} & \multicolumn{1}{c}{2} & \multicolumn{1}{c}{3} & \multicolumn{1}{c}{4} & \multicolumn{1}{c}{5} & \multicolumn{1}{c}{6}\\
\midrule
\textbf{$\left(x, y\right)$} & $\left(\frac{13}{36}, -\frac{13}{36}\right)$ & $\left(\frac{11}{36}, \frac{1}{36}\right)$ &  $\left(\frac{1}{36}, \frac{11}{36}\right)$  & $\left(-\frac{13}{36}, \frac{13}{36}\right)$ & $\left(-\frac{11}{36}, -\frac{1}{36}\right)$ & $\left(-\frac{1}{36}, -\frac{11}{36}\right)$  \\ 
\botrule
\end{tabular}
\caption{$K_{3,3}$ toroidal embedding coordinates}
\label{tab3}
\end{table}

Note that the packing radius for $K_{3,3}$ is slightly less than that of the octahedral graph. This is of course to be expected for two reasons: On one hand, as previously mentionned, the octahedral contact graph is known to be the optimal packing for 6 nodes. On the other hand, $K_{3,3}$ is the 3-regular graph for 6 nodes, whereas the octahedral graph is the 4-regular graph for 6 nodes. This implies that the octahedral graph has greater connectivity than $K_{3,3}$ and therefore also a greater packing radius on the flat torus. We now propose the following conjecture:

\begin{conjecture}
{
  The herein proposed embedding of $K_{3,3}$, with coordinates as described in Table~\ref{tab3} is unique up to an isometry.
}
\end{conjecture}

\section{\texorpdfstring{$K_{6}$}{Lg} and \texorpdfstring{$K_{7}$}{Lg} Toroidal Embeddings}\label{sec8}

\begin{figure}[ht]%
\centering
\includegraphics[width=0.4\textwidth]{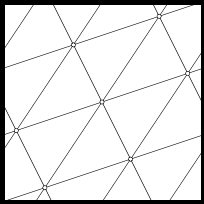}
\caption{\texorpdfstring{$K_{7}$}{Lg}  Toroidal embeddings wihout crossing edges}\label{fig7}
\end{figure}

In this final section, we simply illustrate that the next complete graphs $K_{6}$ and $K_{7}$ can also be embedded on the flat torus, although not as penny graphs. This is shown in Figure~\ref{fig7} with a toroidal embedding of $K_{7}$, as a section of a tilted triangular lattice. Indeed, it makes sense that $K_{7}$ is a section of a triangular lattice because it is the 6-regular graph with 7 nodes, and all nodes of the triangular lattice have degree 6. Removing any single node in Figure~\ref{fig7} would yield a toroidal embedding of $K_{6}$.

\section{Conclusion}\label{sec6}

The exploration of toroidal penny graphs presents exciting possibilities for advancing complex graph visualization techniques. By leveraging the unique properties of penny graphs and their embeddings on the flat torus, researchers can develop more effective methods for representing intricate relationships within data. As our understanding of toroidal penny graphs deepens, their application in visualizing complex graphs may prove invaluable, enhancing clarity and insight in various fields, including network analysis, data science, and computational geometry. Continued investigation in this area holds the promise of innovative approaches to graph visualization, ultimately contributing to more robust analytical tools and methodologies.

\bibliography{k33penny}


\begin{thebibliography}{16}
\ifx \bisbn   \undefined \def \bisbn  #1{ISBN #1}\fi
\ifx \binits  \undefined \def \binits#1{#1}\fi
\ifx \bauthor  \undefined \def \bauthor#1{#1}\fi
\ifx \batitle  \undefined \def \batitle#1{#1}\fi
\ifx \bjtitle  \undefined \def \bjtitle#1{#1}\fi
\ifx \bvolume  \undefined \def \bvolume#1{\textbf{#1}}\fi
\ifx \byear  \undefined \def \byear#1{#1}\fi
\ifx \bissue  \undefined \def \bissue#1{#1}\fi
\ifx \bfpage  \undefined \def \bfpage#1{#1}\fi
\ifx \blpage  \undefined \def \blpage #1{#1}\fi
\ifx \burl  \undefined \def \burl#1{\textsf{#1}}\fi
\ifx \doiurl  \undefined \def \doiurl#1{\url{https://doi.org/#1}}\fi
\ifx \betal  \undefined \def \betal{\textit{et al.}}\fi
\ifx \binstitute  \undefined \def \binstitute#1{#1}\fi
\ifx \binstitutionaled  \undefined \def \binstitutionaled#1{#1}\fi
\ifx \bctitle  \undefined \def \bctitle#1{#1}\fi
\ifx \beditor  \undefined \def \beditor#1{#1}\fi
\ifx \bpublisher  \undefined \def \bpublisher#1{#1}\fi
\ifx \bbtitle  \undefined \def \bbtitle#1{#1}\fi
\ifx \bedition  \undefined \def \bedition#1{#1}\fi
\ifx \bseriesno  \undefined \def \bseriesno#1{#1}\fi
\ifx \blocation  \undefined \def \blocation#1{#1}\fi
\ifx \bsertitle  \undefined \def \bsertitle#1{#1}\fi
\ifx \bsnm \undefined \def \bsnm#1{#1}\fi
\ifx \bsuffix \undefined \def \bsuffix#1{#1}\fi
\ifx \bparticle \undefined \def \bparticle#1{#1}\fi
\ifx \barticle \undefined \def \barticle#1{#1}\fi
\bibcommenthead
\ifx \bconfdate \undefined \def \bconfdate #1{#1}\fi
\ifx \botherref \undefined \def \botherref #1{#1}\fi
\ifx \url \undefined \def \url#1{\textsf{#1}}\fi
\ifx \bchapter \undefined \def \bchapter#1{#1}\fi
\ifx \bbook \undefined \def \bbook#1{#1}\fi
\ifx \bcomment \undefined \def \bcomment#1{#1}\fi
\ifx \oauthor \undefined \def \oauthor#1{#1}\fi
\ifx \citeauthoryear \undefined \def \citeauthoryear#1{#1}\fi
\ifx \endbibitem  \undefined \def \endbibitem {}\fi
\ifx \bconflocation  \undefined \def \bconflocation#1{#1}\fi
\ifx \arxivurl  \undefined \def \arxivurl#1{\textsf{#1}}\fi
\csname PreBibitemsHook\endcsname

\bibitem[\protect\citeauthoryear{Kuratowski}{1930}]{bib0}
\begin{barticle}
\bauthor{\bsnm{Kuratowski}, \binits{C.}}:
\batitle{Sur le problème des courbes gauches en topologie}.
\bjtitle{Fundamenta Mathematicae}
\bvolume{15}(\bissue{1}),
\bfpage{271}--\blpage{283}
(\byear{1930})
\end{barticle}
\endbibitem

\bibitem[\protect\citeauthoryear{Hopcroft and Tarjan}{1974}]{bib1}
\begin{barticle}
\bauthor{\bsnm{Hopcroft}, \binits{J.}},
\bauthor{\bsnm{Tarjan}, \binits{R.}}:
\batitle{Efficient planarity testing}.
\bjtitle{J. ACM}
\bvolume{21}(\bissue{4}),
\bfpage{549}--\blpage{568}
(\byear{1974})
\doiurl{10.1145/321850.321852}
\end{barticle}
\endbibitem

\bibitem[\protect\citeauthoryear{Harborth}{1974}]{bib2}
\begin{barticle}
\bauthor{\bsnm{Harborth}, \binits{H.}}:
\batitle{Lösung zu problem 664a}.
\bjtitle{Elem. Math.}
\bvolume{29},
\bfpage{14}--\blpage{15}
(\byear{1974})
\end{barticle}
\endbibitem

\bibitem[\protect\citeauthoryear{Eades and Wormald}{1990}]{bib3}
\begin{barticle}
\bauthor{\bsnm{Eades}, \binits{P.}},
\bauthor{\bsnm{Wormald}, \binits{N.C.}}:
\batitle{Fixed edge-length graph drawing is np-hard}.
\bjtitle{Discrete Applied Mathematics}
\bvolume{28}(\bissue{2}),
\bfpage{111}--\blpage{134}
(\byear{1990})
\doiurl{10.1016/0166-218X(90)90110-X}
\end{barticle}
\endbibitem

\bibitem[\protect\citeauthoryear{Eades and Whitesides}{1996}]{bib4}
\begin{barticle}
\bauthor{\bsnm{Eades}, \binits{P.}},
\bauthor{\bsnm{Whitesides}, \binits{S.}}:
\batitle{The logic engine and the realization problem for nearest neighbor graphs}.
\bjtitle{Theoretical Computer Science}
\bvolume{169}(\bissue{1}),
\bfpage{23}--\blpage{37}
(\byear{1996})
\doiurl{10.1016/S0304-3975(97)84223-5}
\end{barticle}
\endbibitem

\bibitem[\protect\citeauthoryear{Tammes}{1930}]{bib9}
\begin{bbook}
\bauthor{\bsnm{Tammes}, \binits{P.M.L.}}:
\bbtitle{On the Origin of Number and Arrangement of the Places of Exit on the Surface of Pollen-grains.}
\bpublisher{J.H. De Bussy.},
\blocation{NY}
(\byear{1930})
\end{bbook}
\endbibitem

\bibitem[\protect\citeauthoryear{Leech.}{1957}]{bib10}
\begin{barticle}
\bauthor{\bsnm{Leech.}, \binits{J.}}:
\batitle{Equilibrium of sets of particles on a sphere.}
\bjtitle{The Mathematical Gazette}
\bvolume{41}(\bissue{336}),
\bfpage{81}--\blpage{90}
(\byear{1957})
\end{barticle}
\endbibitem

\bibitem[\protect\citeauthoryear{Kottwitz}{1991}]{bib12}
\begin{barticle}
\bauthor{\bsnm{Kottwitz}, \binits{D.A.}}:
\batitle{The densest packing of equal circles on a sphere}.
\bjtitle{Acta Crystallographica Section A}
\bvolume{47},
\bfpage{158}--\blpage{165}
(\byear{1991})
\end{barticle}
\endbibitem

\bibitem[\protect\citeauthoryear{Clare and Kepert}{1986}]{bib11}
\begin{barticle}
\bauthor{\bsnm{Clare}, \binits{B.W.}},
\bauthor{\bsnm{Kepert}, \binits{D.L.}}:
\batitle{The closest packing of equal circles on a sphere}.
\bjtitle{Proceedings of the Royal Society of London. A. Mathematical and Physical Sciences}
\bvolume{405},
\bfpage{329}--\blpage{344}
(\byear{1986})
\end{barticle}
\endbibitem

\bibitem[\protect\citeauthoryear{Ballinger et~al.}{2006}]{bib13}
\begin{botherref}
\oauthor{\bsnm{Ballinger}, \binits{B.}},
\oauthor{\bsnm{Blekherman}, \binits{G.}},
\oauthor{\bsnm{Cohn}, \binits{H.}},
\oauthor{\bsnm{Giansiracusa}, \binits{N.}},
\oauthor{\bsnm{Kelly}, \binits{E.}},
\oauthor{\bsnm{Schuermann}, \binits{A.}}:
Experimental study of energy-minimizing point configurations on spheres.
Experimental Mathematics
\textbf{18}
(2006)
\doiurl{10.1080/10586458.2009.10129052}
\end{botherref}
\endbibitem

\bibitem[\protect\citeauthoryear{Bachoc and Vallentin}{2009}]{bib14}
\begin{barticle}
\bauthor{\bsnm{Bachoc}, \binits{C.}},
\bauthor{\bsnm{Vallentin}, \binits{F.}}:
\batitle{Optimality and uniqueness of the (4,10,1/6) spherical code}.
\bjtitle{Journal of Combinatorial Theory, Series A}
\bvolume{116},
\bfpage{195}--\blpage{204}
(\byear{2009})
\doiurl{10.1016/j.jcta.2008.05.001}
\end{barticle}
\endbibitem

\bibitem[\protect\citeauthoryear{Sloane et~al.}{1998}]{bib15}
\begin{botherref}
\oauthor{\bsnm{Sloane}},
\oauthor{\bsnm{Hardin}},
\oauthor{\bsnm{Vallentin}},
\oauthor{\bsnm{Smith}}, et al.:
Tables of spherical codes.
published electronically at, NeilSloane.com/packings/
(1998)
\end{botherref}
\endbibitem

\bibitem[\protect\citeauthoryear{Leech and Sloane}{1999}]{bib16}
\begin{bbook}
\bauthor{\bsnm{Leech}, \binits{J.}},
\bauthor{\bsnm{Sloane}, \binits{N.J.A.}}:
\bbtitle{Sphere Packing and Error-Correcting Codes},
pp. \bfpage{136}--\blpage{156}.
\bpublisher{Springer},
\blocation{New York, NY}
(\byear{1999}).
\doiurl{10.1007/978-1-4757-6568-7_5} .
\burl{https://doi.org/10.1007/978-1-4757-6568-7_5}
\end{bbook}
\endbibitem

\bibitem[\protect\citeauthoryear{Dickinson et~al.}{2011}]{bib8}
\begin{barticle}
\bauthor{\bsnm{Dickinson}, \binits{W.}},
\bauthor{\bsnm{Guillot}, \binits{D.}},
\bauthor{\bsnm{Keaton}, \binits{A.}},
\bauthor{\bsnm{Xhumari}, \binits{S.}}:
\batitle{Optimal packings of up to five equal circles on a square flat torus}.
\bjtitle{Beiträge zur Algebra und Geometrie / Contributions to Algebra and Geometry}
\bvolume{52},
\bfpage{315}--\blpage{333}
(\byear{2011})
\doiurl{10.1007/s13366-011-0029-7}
\end{barticle}
\endbibitem

\bibitem[\protect\citeauthoryear{Musin and Nikitenko}{2016}]{bib6}
\begin{botherref}
\oauthor{\bsnm{Musin}, \binits{O.}},
\oauthor{\bsnm{Nikitenko}, \binits{A.}}:
Optimal packings of congruent circles on a square flat torus.
Discrete \& Computational Geometry
\textbf{55}
(2016)
\doiurl{10.1007/s00454-015-9742-6}
\end{botherref}
\endbibitem

\bibitem[\protect\citeauthoryear{Connelly et~al.}{2016}]{bib7}
\begin{botherref}
\oauthor{\bsnm{Connelly}, \binits{R.}},
\oauthor{\bsnm{Funkhouser}, \binits{M.}},
\oauthor{\bsnm{Kuperberg}, \binits{V.}},
\oauthor{\bsnm{Solomonides}, \binits{E.}}:
Packings of equal disks in a square torus
(2016).
\url{https://arxiv.org/abs/1512.08762}
\end{botherref}
\endbibitem

\end{thebibliography}

\end{document}